\documentclass[12pt,indentfirst]{article}

\usepackage{jsr}

\def\cl#1{\begin{center}\Large #1\end{center}}
\def\ce#1{{\centering #1 \par}}
\def\Xt{\widetilde{X}}
\def\Xh{\widehat{X}}
\def\Ht{\widetilde{H}}
\def\Mt{\widetilde{M}}
\def\Tt{\widetilde{T}}
\def\Rt{\widetilde{R}}
\def\Bproc{$\{B_t\}$}
\def\Xproc{$\{X_t\}$}
\def\Xtproc{$\{\Xt_t\}$}
\def\Xhproc{$\{\Xh_t\}$}
\def\Rproc{$\{R_t\}$}
\def\Rtproc{$\{\Rt_t\}$}
\def\Geometric{{\rm Geometric}}
\def\pibar{\overline{\pi}}

\begin{document}

\baselineskip=18pt

\cl{Hitting Time and Convergence Rate Bounds \\
for Symmetric Langevin Diffusions}

\ce{by}

\ce{\gareth\ and \jeff}

\medskip \ce{(November, 2016)}

\bigskip

\abstract{We provide quantitative bounds on the convergence to
stationarity of
real-valued Langevin diffusions with symmetric target densities.}

\newpartitle{Keywords} Langevin diffusion; computable bounds;
coupling; stochastic monotonicity.

\section{Introduction}

Quantitative (computable) bounds on the convergence of Markov processes
to stationarity are an important and widely studied topic, particularly
in the context of Markov chain Monte Carlo (MCMC) algorithms (see e.g.\
Roberts and Tweedie, 1999; Rosenthal, 1995b, 1996, 2002; Jones and
Hobert, 2001, 2004; Baxendale, 2005; and references therein).  Most of
this research has focused on discrete-time Markov chains.  However,
continuous-time Langevin diffusions also converge to stationary
distributions, and quantitative bounds on their convergence are also
of interest.  A start in this direction was made for a few specific
examples by Roberts and Rosenthal (1996) and Roberts and Tweedie
(2000), but much remains to be done.

The current paper was inspired by a question from John Lafferty
(personal communication), who asked about quantitative convergence
bounds for Langevin diffusions on $\IR$ with target densities
proportional to $e^{-|x|^\beta}$ for some fixed $\beta > 1$.  Below we
provide quantitative convergence upper bounds for such diffusions,
and more generally for any symmetric real Langevin diffusion satisfying
certain conditions.  Our bounds are conservative, but are still numerically
modest.  For example, we show that in the $e^{-|x|^\beta}$ case with
$\beta=2$, if we begin the diffusion at $y=2$, then it converges to
within 0.01 of stationarity in total variation distance by time 20.
(Or, if $\beta=1.1$ and $y=10$, then it is within 0.01 of stationarity
by time 34.)

Our proof requires bounds on hitting times, which are also developed
below using probability generating functions.  It uses the coupling
inequality, and the stochastic monotonicity of the diffusions,
following the general approach of Lund et al.\ (1996a, 1996b).

\section{Assumptions and Hitting Probabilities}

\def\P{{\bf P}}
\def\E{{\bf E}}
\def\IR{{\bf R}}
\def\L{{\cal L}}
\def\half{{1 \over 2}}
\def\grad{\nabla}
\def\un{\underbar}

Let $\pi:\IR\to[0,\infty)$ be a target density on $\IR$, satisfying
the following:

\medskip \itemitem{\bf (A1)}
(i) $\pi$ is {\it symmetric}, i.e.\ $\pi(-x)=\pi(x)$ for all $x\in\IR$;
\hfil\break
(ii) $\pi$ is $C^1$, i.e.\ is continuously differentiable;
\hfil\break
(iii) $-\grad \log \pi(x) \ge 0$ for all $x \ge 0$;
\hfil\break
(iv) there is $b>0$ such that $-\grad \log \pi(x) \ge b >
0$ for all $x \ge 1$.

\medskip\noindent
For example, (A1) is satisfied if $\pi(x) \propto e^{-|x|^\beta}$ for any
fixed $\beta > 1$,
in which case for $x \ge 0$ we have
$- \grad \log \pi(x) = - \grad (-x^\beta)
= \beta x^{\beta-1}$, which is $\ge 0$ for $x \ge 0$,
and which is non-decreasing on $[1,\infty)$ so
we can take $b = - \grad \log \pi (1) = \beta > 0$.
\medskip

Let \Xproc\ be a Langevin diffusion for $\pi$, so $dX_t = \half
\, \grad \log \pi(X_t) \, dt + dB_t$ (where $\{B_t\}$ is standard
Brownian motion).  Let $H_y$ be this diffusion's first hitting time
of 0, conditional on starting at $X_0=y$.  We wish to bound the tail
probabilities of $H_y$.

Our key computation is the following bound on the probability
generating function of $H_y$, i.e.\ of $M_y(s) := \E(s^{H_y})$.
We require that:

\medskip \itemitem{\bf (A2)}
The value $s\in\IR$ satisfies that
$s>1$, and $s < \exp(b^2/2)$, and
$$
1 < { \exp\left( b - \sqrt{b^2 - 2\log s} \right)
\over \cos(\sqrt{2\log s}) } < 2
\, .
$$
\medskip

\pnew{Lemma} \plabel{momgenlem}
If $\pi$ satisfies (A1), then for any $y \in \IR$,
and any $s$ satisfying (A2),
the probability generating function $M_y(s) := \E(s^{H_y})$ of the
Langevin diffusion for $\pi$ satisfies
$M_y(s) \le B(\max(1,|y|),s,b)$,
where
$$
B(y,s,b) \ = \
{ \exp\left( (y-1) \Big[ b - \sqrt{b^2 - 2\log s} \Big] \right)
\over \cos(\sqrt{2\log s}) }
\Biggm/ \ \left[ 2 -
{\exp\left(b - \sqrt{b - 2\log s}\right)
\over \cos(\sqrt{2\log s})} \right]
\, .
$$

\bigskip

Lemma~\ref{momgenlem} is proved in Section~\ref{sec-lemmaproof} below.


Assuming Lemma~\ref{momgenlem}, we immediately obtain a bound on the
tail probabilities of the hitting time of 0, as follows.

\pnew{Proposition}
\plabel{hitprop}
If $\pi$ satisfies (A1), then
for any $y \in \IR$ and $t>0$,
and any $s$ satisfying (A2),
the hitting time $H_y$ of the Langevin diffusion for $\pi$ satisfies
$$
\P(H_y \ge t) \ \le \ s^{-t} \, B\Big( \max(1,|y|), \, s \Big)
\, .
$$

\proof It follows by Markov's inequality that
$$
\P(H_{y} \ge t)
\ = \ \P(s^{H_y} \ge s^t)
\ \le \ s^{-t} \, \E(s^{H_y})
\ = \ s^{-t} \, M_{y}(s)
\ \le \ s^{-t} \, B(\max(1,|y|),s,b)
\, .
\eqqed
$$


\medskip
\newpartitle{Numerical Example}
Suppose $\pi(x) \propto e^{-|x|^\beta}$, so $b=\beta$ as above.
Then if, say, $y=\beta=b=2$,
then choosing $s=1.4$, we compute numerically
from Proposition~\ref{hitprop} that
$\P(H_{y} \ge t) <0.01$ whenever $t \ge 20$.
This indicates that this
process has probability over 99\% of hitting 0 by time 20.
By contrast, if $\beta=1.1$ and $y=10$, then taking $s=1.3$ gives that
$\P(H_{y} \ge t) <0.01$ whenever $t \ge 34$.

\section{Convergence Bounds: Reflected Case}

Our main interest is in bounds on the convergence time to stationarity.
We achieve this via the coupling inequality, using {\it stochastic
monotonicity} as in Lund et al.\ (1996a, 1996b).

We first consider the version of the process which is
``reflected'' or ``folded'' at zero, which is equivalent to
considering just the {\it absolute value} of the process.
That is, we consider the process $R_t := |X_t|$.  This process behaves
just like \Xproc\ on the positive half-line.  But when it hits 0, it
reflects back to the positive part rather than go negative.

By symmetry, this process has a stationary density $\pibar$ equal to
$\pi$ restricted to the positive half-line, i.e.\ $\pibar(x) = 2 \,
\pi(x) \one_{x \ge 0}$.
We seek specific quantitative computable bounds on the total variation
distance to stationarity of this process after time $t$, i.e.\ on
$$
\|\L_y(R_t) - \pibar\|
\ := \ \sup_{A \subseteq \IR} |\P_y(R_t \in A) - \pibar(A)|
\, ,
$$
where the supremum is taken over all measurable subsets $A \subseteq
\IR$, and the subscript $y$ indicates that the process was
started at the state $y$, and $\pibar(A) := \int_A \pibar(x) \, dx$.
We shall prove the following.

\pnew{Proposition}
\plabel{refprop}
If $\pi$ satisfies (A1), then for
the reflected Langevin diffusion for $\pi$, started at state $y \ge 1$, 
for any $s$ satisfying (A2),
the total variation distance to stationarity at time $t>0$ satisfies
$$
\|\L_y(R_t) - \pibar\|
\ \le \
2 \, \pi[0,y] \, s^{-t} \, B(y,s,b)
\, + 2 \, \int_{y}^\infty \pi(z) \, s^{-t} \, B(z,s,b) \, dz
\, ,
$$
with $B(y,s,b)$ as in Lemma~\ref{momgenlem}.

\proof
Suppose we begin two separate copies of the process,
\Rproc\ and \Rtproc, started at $y$ and at $z$, respectively,
and couple them by using the same Brownian
motion \Bproc\ for both of them.  Then the processes'
relative ordering is preserved.
It follows that if the {\it larger} of the two processes hits~0,
then the other process must also equal~0 at that same time.  Hence,
the two processes must couple (i.e., become equal) by the time the
larger of the two processes hits~0.

Therefore, if the coupling time is
$U$, then $U$ is stochastically bounded above by the hitting time
$H_{\max(y,z)}$ of 0 from the state $\max(y,z)$, i.e.\
$\P(U \ge t) \le \P(H_{\max(y,z)} \ge t)$.  Hence, by the usual
coupling inequality (e.g.\ Roberts and Rosenthal, 2004, Section~4.1),
the total variation distance of the processes after time $t$ satisfies
$$
\|\L_y(R_t) - \L_z(\Rt_t)\|
\ := \ \sup_{A \subseteq \IR} |\P(R_t \in A) - \P(\Rt_t \in A)|
$$
$$
\ \le \ \P(U \ge t)
\ \le \ \P(H_{\max(y,z)} \ge t)
\, .
$$
Therefore, by Proposition~\ref{hitprop}, for any $s \ge 1$,
$\|\L_y(R_t) - \L_z(\Rt_t)\| \, \le \, s^{-t} \, B(\max(1,y,z),s,b)$.

Suppose now that we start the \Rproc\ process at $y \ge 1$,
and start the \Rtproc\ process at a state $Z \sim \pibar$ chosen
randomly from the stationary distribution.  Then
$$
\|\L_y(R_t)-\pibar\|
\ \equiv \ \Big\| \L_y(R_t) - \Big( \E_{Z \sim \pibar} \L_Z(\Rt_t) \Big) \Big\|
$$
$$
\ \le \ \E_{Z \sim \pibar} \|\L_y(R_t) - \L_Z(\Rt_t)\|
$$
$$
\ \le \ \E_{Z \sim \pibar} [\P(H_{\max(y,Z)} \ge t)]
$$
$$
\ \le \ \E_{Z \sim \pibar} [s^{-t} \, B(\max(1,y,Z),s,b)]
$$
$$
\ = \ \int_0^\infty \pibar(z) \, s^{-t} \, B(\max(1,y,z),s,b) \, dz
$$
$$
\ = \
\pibar[0,y] \, s^{-t} \, B(y,s,b)
\, + \, \int_{y}^\infty \pibar(z) \, s^{-t} \, B(z,s,b) \, dz
$$
$$
\ = \
2 \, \pi[0,y] \, s^{-t} \, B(y,s,b)
\, + 2 \, \int_{y}^\infty \pi(z) \, s^{-t} \, B(z,s,b) \, dz
\, .
\eqqed
$$

\newpartitle{Numerical Example}
Again let $\pi(x) \propto e^{-|x|^\beta}$.
Then from the above,
$$
\|\L_y(R_t)-\pibar\| \ \le \
2 \, \pi[0,y] \, s^{-t} \, B(y,s,b)
\, + 2 \, \int_{y}^\infty \pi(z) \, s^{-t} \, B(z,s,b) \, dz
\, .
$$
Suppose again that $y=\beta=b=2$.  Then choosing $s=1.4$,
Proposition~\ref{refprop} shows that
$\|\L_y(R_t)-\pibar\|<0.01$ whenever $t \ge 20$, i.e.\ the process also
converges to within 99\% of its stationarity distribution by time 20.
By contrast, with $y=10$ and $\beta=b=1.1$, we find choosing $s=1.3$ that
$\|\L_y(R_t)-\pibar\|<0.01$ whenever $t \ge 34$.

\section{Convergence Bounds: Unreflected Case}

Finally, we consider convergence bounds on the full, unreflected
diffusion \Xproc.  Here we cannot use stochastic monotonicity
directly, because the full diffusion has no lowest state (or even
lower bound) on which to force two copies of the diffusion to couple.


Nevertheless, using the symmetry condition A1(i), we are able to prove
that the same convergence time bounds hold for the unreflected case
as for the reflected case:

\pnew{Proposition}
\plabel{unrefprop}
If $\pi$ satisfies (A1), then for
the full unreflected Langevin diffusion for $\pi$, started at state $y \ge 1$, 
for any $s$ satisfying (A2),
the total variation distance to stationarity at time $t>0$ satisfies
$$
\|\L_y(X_t) - \pi\|
\ \le \
2 \, \pi[0,y] \, s^{-t} \, B(y,s,b)
\, + 2 \, \int_{y}^\infty \pi(z) \, s^{-t} \, B(z,s,b) \, dz
\, ,
$$
with $B(y,s,b)$ as in Lemma~\ref{momgenlem}.

\proof We jointly define two copies \Xproc\ and \Xtproc\ of the
Langevin diffusion for $\pi$, by $X_0=y$, $\Xt_0 \sim \pi$,
$dX_t = \half \, \grad \log \pi(X_t) \, dt + dB_t$,
and $d\Xt_t = \half \, \grad \log \pi(\Xt_t) \, dt - dB_t$,
where $\{B_t\}$ is the same standard Brownian motion in both cases.
In particular, \Xproc\ and \Xtproc\ are {\it anti-coupled}, driven
by the same Brownian motion with opposite signs.

For this joint process, let
$$
\tau \ = \ \inf\{t \ge 0 : X_t = \Xt_t\}
$$
be the first time they meet.  Thus, $\tau$ is a stopping time for the
joint process.  Finally, define another process \Xhproc\ by
$\Xh_0 = \Xt_0$, and
$$
d\Xh_t \ = \
\left\{
\begin{array}{lr}
\half \, \grad \log \pi(\Xt_t) \, dt - dB_t \, , & \quad t \le \tau \\
\half \, \grad \log \pi(\Xt_t) \, dt + dB_t \, , & \quad t > \tau
\end{array}
\right.
$$
That is, \Xhproc\ is the same as \Xtproc\ up to the meeting time
$\tau$, after which \Xhproc\ is the same as \Xproc.  (This
construction is valid since $\tau$ is a joint stopping time.)

We claim that this joint process preserves the absolute-value ordering
of the processes \Xproc\ and \Xhproc.  That is,
if $|X_0| \ge |\Xh_0|$ then $|X_t| \ge |\Xh_t|$ for all $t \ge 0$,
while if $|X_0| \le |\Xh_0|$ then $|X_t| \le |\Xh_t|$ for all $t \le 0$.
Indeed, since the diffusions have continuous sample paths, this is
immediately true without the absolute value signs.  The only remaining
case is if $X_t$ and $\Xh_t$ are of opposite sign.  But in that case,
we must have $t < \tau$, so that \Xproc\ and \Xhproc\ are driven by
Brownian motions of opposite sign.  By the symmetry condition A1(i),
this means that their
absolute values are driven by the same Brownian motion, and hence
again cannot cross because of the continuous sample paths.  So, either
way, the absolute-value ordering is preserved.

The rest of the argument is identical to that of
Proposition~\ref{refprop}.  Indeed, when the larger (in absolute
value) of the two processes reaches zero, the smaller one must also
reach zero, so they must have coupled by that time.  That is,
conditional on $X_0=y$ and $\Xt_0=z$, we must have $\tau \le H_{\max(y,z)}$.

Hence, if we start the \Xproc\ process at a state $y \ge 1$,
and start the \Xtproc\ process at a state $Z \sim \pi$ chosen
randomly from the stationary distribution, then
applying Proposition~\ref{hitprop}, we must again have
$$
\|\L_y(X_t)-\pi\|
\ \equiv \ \Big\| \L_y(X_t) - \Big( \E_{Z \sim \pi} \L_Z(\Xt_t) \Big) \Big\|
$$
$$
\ \le \ \E_{Z \sim \pi} \|\L_y(X_t) - \L_Z(\Xt_t)\|
$$
$$
\ \le \ \E_{Z \sim \pi} [\P(H_{\max(|y|,|Z|)} \ge t)]
$$
$$
\ = \ \E_{Z \sim \pibar} [\P(H_{\max(|y|,Z)} \ge t)]
$$
$$
\ \le \
2 \, \pi[0,y] \, s^{-t} \, B(y,s,b)
\, + 2 \, \int_{y}^\infty \pi(z) \, s^{-t} \, B(z,s,b) \, dz
\, ,
$$
exactly as in the proof of Proposition~\ref{refprop}.
\qed

\newpartitle{Numerical Example}
Since the bounds in Proposition~\ref{unrefprop}
are the same as those in Proposition~\ref{refprop}, the same bounds
still apply.  For example,
if $y=\beta=b=2$, then choosing $s=1.4$,
Proposition~\ref{unrefprop} again shows that
$\|\L_y(X_t)-\pibar\|<0.01$ whenever $t \ge 20$, i.e.\ the full
unreflected process also
converges to within 99\% of its stationarity distribution by time 20.
Or, with $y=10$ and $\beta=b=1.1$, Proposition~\ref{unrefprop}
with $s=1.3$ yields that
$\|\L_y(X_t)-\pibar\|<0.01$ whenever $t \ge 34$.

\remark It is perhaps surprising that we obtain identical convergence
bounds in both the reflected and the unreflected (original) case.
Indeed, in general the reflected case should converge faster.
However, since our upper bounds are derived using the hitting time
of the larger process (in absolute value) to reach zero, we
obtain the same bound in both cases.

\section{Proof of Lemma~\ref{momgenlem}}
\label{sec-lemmaproof}

Finally, we prove Lemma~\ref{momgenlem}.  We assume $y \ge 1$; the
case $y \le -1$ then follows by symmetry, and the case $|y| < 1$
then follows since, by monotonicity, the hitting time of~0 from such $y$
is stochastically bounded above by the hitting time of~0 from~1.


We first introduce an indicator process $\{I_t\}_{t \ge 0}$, defined
as follows.  We begin with $I_0=0$.  Then, each time $X_t$ hits 1,
we set $I_t=1$.  Then, each time $X_t$ hits 2 or 0, we set $I_t=0$.
Intuitively, $I_t$ indicates whether at time $t$ we are ``drifting
towards 1'', or are ``waiting to hit 0 or 2 from 1''.

Then, we let \Xtproc\ be a slight modification of $\{X_t\}$,
as follows.  \Xtproc\ mostly follows the same dynamics as \Xproc.
However, whenever $I_t=1$, the drift of \Xtproc\ is instead 0, i.e.\
we replace the Langevin diffusion dynamics by standard brownian motion.
Also, when $I_t=0$, the drift of \Xtproc\ is instead the constant
value $-b$.


Now, because of the assumptions (A1),
this new process is {\it stochastically larger} than the original
process up to time $H_y$, i.e.\ it can only take {\it longer} to hit 0.
So, writing $\Ht_y$ for the first time the modified process
\Xtproc\ hits 0, and $\Mt_y(s) := \E(s^{\Ht_y})$ for its probability
generating function, we must have $M_y(s) \le \Mt_y(s)$ for all $s,y
\ge 0$.  That is, we can (and will) use the hitting time of 0
for the modified process, as an upper bound on the he hitting time
of 0 for the original process.
So, it suffices to show that $\Mt_y(s) = B(y,s,b)$, which we now do.

For the modified process $\{\Xt_t\}$, we break up the journey from $y \ge 1$
to 0 into steps:
\hfil\break {\bf 1.} Reach the state 1.
\hfil\break {\bf 2.} From there, reach either state 0 or state 2.
\hfil\break {\bf 3(a).} If it reached state 0, then we're done.
\hfil\break {\bf 3(b).} If instead it reached state 2, then return to step~1.
\par\medskip
Now, since $\{\Xt_t\}$ has zero drift on $[0,2]$, it follows that
after reaching state~1, the process has equal probability 1/2 of
reaching either 0 or 2.  Therefore,
the number of times the process will return to step
1 before finally hitting 0 is
a geometric random variable $G \sim \Geometric(1/2)$ with
$\P[G=k] = 2^{-k-1}$ for $k=0,1,2,\ldots$.  It then follows
that the time for $\{\Xt_t\}$ to reach state 0 from $y$ can be written as
$$
\Ht_y \ = \
T_{y1} + \sum_{i=1}^G (\Tt^{(i)}_{12|2} + T^{(i)}_{21}) + \Tt_{12|0}
\, .
$$
Here $T_{y1}$ is the random time for \Xtproc\ to reach 1 from $y$ with
constant drift $-b$,
and each $T^{(i)}_{21}$ is an {\it independent} random time for \Xproc\
to reach 1 from 2 with drift $-b$,
and each $\Tt^{(i)}_{12|2}$ is an {\it independent} random time for \Xproc\
with drift~0 to reach 2 from 1 with drift~0 but
conditional on reaching 2 before 0,
and $\Tt_{10|0}$ is an {\it independent} random time for standard
Brownian motion to reach 0 from 1 with drift~0 but
conditional on reaching 0 before 2.



Hence, with corresponding notation, the probability generating function
$\Mt_y(s) := \E(s^{\Ht_y})$ is given by
\begin{equation}
\label{Meqn}
\Mt_y(s)
\ = \ M_{y1}(s) \times M_G\Big(\Mt_{12|2}(s) \times M_{21}(s)\Big)
	\times \Mt_{12|0}(s)
\, .
\end{equation}
Now, these various formulas are known.
First, by symmetry, $\Mt_{12|2}(s) = \Mt_{12|0}(s) = M_0^*(s)$ where $M_x^*(s)$
is the probability generating function for the time taken by standard Brownian
motion to reach $\pm 1$ when started at $x$ (where $-1 < x < 1$), and
this is known to be given for $s>1$ by
$$
M^*_x(s)
\ = \ {\cos (x \sqrt{2 \log s}) \over \cos (\sqrt{2 \log s} )}
\, ,
$$
so that
$M^*_0(s) = 1 / \cos (\sqrt{2 \log s} )$.
Also, $G$ has known
probability generating function given by $M_G(r) := \E(r^G)
= (1-1/2) / (1-r/2) = 1/(2-r)$ for $1<r<2$.

Finally, it is known (see e.g.\
pp.~295 and 309 of Borodin and Salminen, 1996;
Proposition 3.3.5
of Etheridge, 2002)
that if $\{W_t\}$ is standard Brownian motion, and
$T_{a,b} = \inf\{t \ge 0 : W_t = a+bt\}$, then
for $\theta>0$ and $a>0$ and $b \ge 0$,
$$
\E\left[ \exp(-\theta T_{a,b}) \right]
\ = \
\exp\left( -a \left[ b + \sqrt{b^2+2\theta} \right] \right)
\, .
$$
Hence, with the identification $\theta:=-\log s$ and $a=y-1$,
this indicates that
$$
M_{y1}(s)
\ = \ \exp\left( (y-1) \Big[ b - \sqrt{b^2 - 2 \log s)} \Big] \right)
\, ,
$$
and $M_{21}(s)$ then follows by setting $y=2$.

Plugging these various formulae into \eqref{Meqn},
it follows by direct algebra (verified using the {\it Mathematica}
symbolic algebra software, Wolfram 1988)
that $\Mt_y(s) = B(y,s,b)$ with $B(y,s,b)$ as stated.
Then, by monotonicity, $M_y(s) \le B(y,s,b)$, giving the result.
\qed

\remark An examination of the proof of Lemma~\ref{momgenlem} indicates
that assumption (A1) is not strictly necessarily, and quantitative
bounds could be obtained by similar methods for Langevin diffusions
for other target densities $\pi$ as well.

\section*{References}
\frenchspacing

P.H. Baxendale (2005), Renewal theory and computable convergence rates
for geometrically ergodic Markov chains.
Ann. Appl. Prob. {\bf 15}, 700--738.

A.N.~Borodin and P.~Salminen (1996),
Handbook of Brownian motion : facts and formulae.
Birkh\"auser Verlag, Basel / Boston.

A.~Etheridge (2002),
A course in Financial calculus.
Cambridge University Press.

G.L. Jones and J.P. Hobert (2001), Honest exploration of intractable
probability distributions via Markov chain Monte Carlo. Statistical
Science {\bf 16}, 312--334.

G.L. Jones and J.P. Hobert (2004), Sufficient burn-in for Gibbs samplers
for a hierarchical random effects model.  Ann. Stat. {\bf 32},
784--817.

R.B. Lund and R.L. Tweedie (1996a), Geometric convergence rates of
stochastically ordered Markov chains.  Math. Oper. Research {\bf
21}, 182--194.

R.B. Lund, S.P. Meyn, and R.L. Tweedie (1996b), Computable exponential
convergence rates for stochastically ordered Markov processes.  Ann.
Appl. Prob. {\bf 6}, 218-237.



G.O.~Roberts and J.S.~Rosenthal (1996),
Quantitative bounds for convergence rates of
continuous time Markov processes.
Elec. J. Prob. {\bf 1(9)}, 1--21.

G.O.~Roberts and J.S.~Rosenthal (2004),
General state space Markov chains and MCMC algorithms.
Prob.\ Surv.\ {\bf 1}, 20--71.

G.O. Roberts and R.L. Tweedie (1999),
Bounds on regeneration times and convergence rates for Markov chains.
Stoch. Proc. Appl. {\bf 80}, 211--229.
See also the corrigendum, Stoch. Proc. Appl. {\bf 91} (2001),
337--338.

G.O. Roberts and R.L. Tweedie (2000),
Rates of convergence of stochastically monotone and continuous time
Markov models.  J. Appl. Prob. {\bf 37}, 359--373.

J.S. Rosenthal (1995b), Minorization conditions and convergence rates
for Markov chain Monte Carlo. J. Amer. Stat. Assoc. {\bf 90},
558--566.

J.S. Rosenthal (1996), Convergence of Gibbs sampler for a model
related to James-Stein estimators.  Stat. and Comput. {\bf 6},
269--275.

J.S. Rosenthal (2002),
Quantitative convergence rates of Markov chains: A simple account.
Elec. Comm. Prob. {\bf 7}, No.~13, 123--128.

S. Wolfram (1988), Mathematica: A system for doing mathematics by
computer.  Addison-Wesley, New York.


\end{document}